\documentclass[11pt,a4paper]{smfart} 
\usepackage{mathsmf}
\usepackage{smfthm}
\usepackage[utf8]{inputenc}
\usepackage[T1]{fontenc}
\usepackage[french]{babel}
\usepackage[plainpages]{hyperref}

\renewcommand\theequation{\arabic{equation}}

\title{Une inégalité de Cheeger pour le spectre de Steklov}
\author{Pierre Jammes}
\address{Univ. Nice Sophia Antipolis, CNRS,  LJAD, UMR 7351, 
06100 Nice, France}
\email{pjammes@unice.fr}

\NumberTheoremsIn{}

\begin{document}
\begin{abstract} 
On montre une inégalité de Cheeger pour la première valeur 
propre de Steklov. Elle fait intervenir deux constantes isopérimétriques.
\end{abstract}


\begin{altabstract}
 We prove a Cheeger inequality for the first positive
Steklov eigenvalue. It involves two isoperimetric constants.
\end{altabstract}
\altkeywords{Cheeger inequality, Steklov eigenvalues.}


\maketitle

Soit $M$ une variété compacte à bord et $\gamma\in C^{1,1}(M)$, 
$\rho\in C^0(\partial M)$ deux fonctions densités strictement positives 
sur $M$ et $\partial M$ respectivement.
Le problème aux valeurs propres de Steklov consiste à résoudre l'équation,
d'inconnues $\sigma\in\R$ et $f:\overline M\to\R$,
\begin{equation}
\left\{\begin{array}{ll}
\divergence(\gamma\nabla f)=0 & \textrm{dans }M\\
\gamma\frac{\partial f}{\partial \nu}=\sigma \rho f & \textrm{sur }\partial M
\end{array}\right.
\end{equation}
où $\nu$ est un vecteur unitaire sortant normal au bord. On parle de 
problème de Steklov homogène quand $\gamma\equiv1$ et $\rho\equiv1$, et le
cas $\gamma\not\equiv1$ se rattache au problème de Calder\'on. 
L'ensemble des
réels $\sigma$ solutions du problème forme un spectre discret positif
noté
\begin{equation}
0=\sigma_0(M,g,\rho,\gamma)<\sigma_1(M,g,\rho,\gamma)\leq
\sigma_2(M,g,\rho,\gamma)\ldots
\end{equation}
C'est le spectre d'un opérateur Dirichlet-Neumann $H^1(\partial M)
\to L^2(\partial M)$ défini par $\Lambda_{\rho,\gamma} u=\frac\gamma\rho
\frac{\partial\mathcal H_\gamma u}{\partial\nu}$, où $\mathcal H_\gamma u$ est 
le prolongement harmonique de $u$ pour la densité $\gamma$, c'est-à-dire que 
$\divergence(\gamma\nabla(\mathcal H_\gamma u))=0$. Il est auto-adjoint 
pour la norme de Hilbert $\|u\|^2=\int_{\partial M}u^2\rho\,\de v_g$ (voir
\cite{ba80}, \cite{su90} et \cite{uh09}). Ce spectre possède une 
caractérisation 
variationnelle ; on a en particulier pour la première valeur propre :
\begin{equation}\label{minmax}
\sigma_1(M,g,\rho,\gamma)=\inf_{\stackrel{f\in H^1(\overline M)}%
{\int_{\partial M} f^2\rho\,\de v_g=0}}
\frac{\int_M|\de f|^2\,\gamma\, \de v_g}%
{\int_{\partial M}f^2\,\rho\,\de v_g}.
\end{equation}

L'objet de cet article est de donner une minoration de $\sigma_1(M)$
en fonction d'invariants isopérimétriques qui est analogue à l'inégalité 
de Cheeger pour la première valeur propre du laplacien. Pour ce faire,
on introduit les notations suivantes : si~$D$ est un domaine de $M$,
on note $\partial_E D=\partial D\cap \partial M$ le bord «~extérieur~» de $D$
et $\partial_I D=\partial D\backslash\partial_E D$ son bord intérieur.
Les volumes $n$ et $(n-1)$-dimensionnels seront calculés relativement
aux densités $\gamma$ et $\rho$ et notés $|\cdot|_\gamma$ et 
$|\cdot|_\rho$. On pose alors
\begin{equation}
h(M,g,\gamma)=\inf_{|D|_\gamma\leq\frac{|M|_\gamma}2}
\frac{|\partial_ID|_\gamma}{|D|_\gamma} \textrm{ et }
h'(M,g,\rho,\gamma)=\inf_{|D|_\gamma\leq\frac{|M|_\gamma}2}
\frac{|\partial_ID|_\gamma}{|\partial_ED|_\rho}.
\end{equation}
Quand $\gamma\equiv1$, la constante $h$ est la constante de Cheeger 
classique ; la première valeur propre non nulle du laplacien de Neumann sur 
$M$ est minorée par $h^2/4$ (\cite{ch70},\cite{bu80}). Des analogues dans 
le cas $\gamma\not\equiv1$ ont été introduits dans \cite{br85} et \cite{co97}.
La constante $h'$ intervient spécifiquement dans le problème de Steklov.
Après avoir écrit la première version de cet article, j'ai appris que
J.~Escobar définit une constante presque identique à $h'$ 
dans \cite{es97} dans le but de minorer $\sigma_1$. Mais l'inégalité qui
suit, plus simple et plus explicite que celle de \cite{es97}, semble
avoir échappé aux recherches menées jusqu'à présent :

\begin{theo}
Soit $(M,g)$ une variété riemannienne compacte à bord munie des densités 
$\rho$ et $\gamma$ sur $\partial M$ et $M$. On a $h'(M)>0$ et
$\displaystyle\sigma_1(M)\geq\frac{h(M)\cdot h'(M)}4$.
\end{theo}
\begin{rema}\label{rem1}
On peut en fait montrer deux inégalités différentes : si on modifie 
légèrement les définitions de $h$ et $h'$ en remplaçant
la condition $|D|_\gamma\leq|M|_\gamma/2$ par 
$|\partial_E D|_\rho\leq|\partial M|_\rho/2$, la démonstration
reste valide mais les constantes isopérimétriques ne sont plus les mêmes.
En particulier, la constante $h'$ est remplacée par la constante 
définie par Escobar.
\end{rema}
\begin{rema}
Si on multiplie la métrique par $\lambda^2$ ($\lambda>0$), alors
$\sigma_1$ et $h$ sont divisées par $\lambda$ et $h'$ est invariant. 
Pour des raisons d'homogénéité, on ne peut donc pas avoir une minoration
de $\sigma_1$ de la forme $c\cdot h^\alpha\cdot h'^\beta$ où $c>0$, 
$\alpha\neq1$ et $\beta\in\R$ sont des constantes universelles. Le cas 
$\alpha=1$ et $\beta<1$ est exclu par l'exemple qui suit.
\end{rema}
\begin{exem}\label{ex1}
Soit $(M,g)$ une variété close. Le spectre de Steklov (homogène) de 
$M_n=M\times
[0,\frac1n]$ a été calculé dans \cite{cesg11} (lemme~6.1) et 
$\sigma_1(M_n)\sim\frac1n$ quand $n\to\infty$. On peut encadrer $h(M_n)$ 
uniformément : étant donné un domaine $D$ de $M_n$, on considère $2n$ 
copies de $M_n$ qu'on recolle le long de leurs bords pour obtenir une variété
$M\times S^1$ dont la métrique est indépendante de $n$, et on construit un 
domaine $D'\subset M\times S^1$ en recollant les domaines $D$ par réflexion 
le long du bord. On obtient alors $|D'|=2n|D|$ et $|\partial D'|=
2n|\partial_ID|$ donc
$|\partial_ID|/|D|=|\partial D'|/|D'|\geq h(M\times S^1)$ et $h(M_n)\geq
h(M\times S^1)$. En considérant un domaine de la forme $D=U\times[0,\frac1n]$ 
avec $U\subset M$, on voit que $h(M_n)$ est majoré par une constante.

On peut aussi estimer $h'$: l'exemple d'un domaine $U\times[0,\frac1n]$ 
montre que $h'(M_n)\leq\frac cn$ pour une 
constante $c>0$. Réciproquement, pour un domaine $D$ fixé, $|\partial D_E|$
est indépendant de $n$ et le volume de $\partial D_I$ sur $M_n$ vérifie
$|\partial D_I|_n\geq|\partial D_I|_1/n$. Par conséquent,
$h'(M_n)\geq h'(M_1)/n$. 

On voit que sur cet exemple, l'exposant de $h'$ dans l'inégalité est
optimal.
\end{exem}
\begin{exem}
On considère encore la variété $M\times[0,L]$ où $M$ est une variété
close, mais on suppose que $L\to+\infty$. Dans ce cas, on a $\sigma_1=
1/L$. Parmi les domaines $D$ de volume fixé suffisamment
grand, l'aire de $\partial_ID$ est minimisée par les hypersurfaces
$\partial_ID=M\times\{x\}$, par conséquent $h=2/L$ et $h'=1$. Dans cet
exemple, c'est $h$ qui tend vers~0 quand $\sigma_1\to0$ et pas $h'$.
\end{exem}
\begin{exem}\label{ex2}
Dans \cite{gp10} (section~2), A.~Girouard et I.~Polterovich construisent 
une famille de domaines du plan (formés de deux disques reliés par une anse 
fine) pour lesquels $\sigma_1$, $h$ et $h'$
tendent vers zéro. Contrairement à la remarque et aux exemples qui précèdent, 
le volume reste borné inférieurement et supérieurement.
\end{exem}
\begin{exem}
Les «~haltères de Cheeger généralisées~» construites par P.~Guérini dans 
\cite{gu04} fournissent des exemples à mi-chemin entre les exemples~\ref{ex1} 
et~\ref{ex2}; on va ici étudier le cas de la dimension~3. Étant donné
un $\varepsilon>0$, dans
$\R^3$ muni d'un repère $(Oxyz)$ on considère le cercle unité $S^1$ et le
disque unité $D^2$ du plan $(Oxy)$, et on définit les domaines 
$V_1=\{p\in\R^3,\ d(p,S^1)<\frac14\}$ et $V_2=D^2\times
]-\varepsilon,\varepsilon[$. La variété définie par $M_\varepsilon=V_1\cup V_2$ 
est homéomorphe à une boule et sa constante isopérimétrique $h'$ tend vers~0
quand $\varepsilon\to0$ (considérer un domaine $D\subset M_\varepsilon$
de la forme $U\times]-\varepsilon,\varepsilon[$ avec $U\subset (Oxy)$ contenu
dans le disque de centre $O$ et rayon $3/4$). On a aussi 
$\sigma_1(M_\varepsilon)\to0$ : il suffit de fixer une fonction test 
$(Oxy)\to\R$ d'intégrale nulle et à support dans le disque de centre $O$ et 
rayon $3/4$ et de la relever à $\R^3$. Dans la formule~(\ref{minmax}), 
le numérateur sera proportionnel à $\varepsilon$ tandis que le 
numérateur reste constant. En construisant de la même manière un espace
test de dimension arbitrairement grande, on montre qu'en fait toutes
les valeurs propres tendent vers~0. Ce phénomène d'«~effondrement du spectre~»
avait déjà été observé par A.~Girouard et I.~Polterovich sur l'exemple
précédent.

\end{exem}
\begin{proof}[Démonstration du théorème]
On fixe une métrique $g$ et des densités $\rho$ et $\gamma$. Pour montrer 
que $h'>0$ il
suffit de traiter le cas homogène, le cas général s'y ramène car les
densités sont continues, donc encadrées par des constantes strictement
positives. On se donne un $\varepsilon\in]0,1]$ dont la valeur précise 
sera fixée
plus tard. Si $|\partial_ED|\leq(1-\varepsilon)|\partial M|$, où
$D$ est un domaine intervenant dans la définition de $h'$,
J.~Escobar montre dans \cite{es99} que le rapport $|\partial_ID|/|\partial_ED|$
est minoré par une constante $K(\varepsilon)>0$ (il ne traite que le cas
$\varepsilon=\frac12$ mais sa démonstration se généralise immédiatement).
Si $|\partial_ED|\geq(1-\varepsilon)|\partial M|$ on utilise le fait que
le volume de $M\backslash D$ est minoré par $|M|/2$ pour tout $D$, et donc
qu'il existe une constante $C(g)>0$ telle que $|\partial(M\backslash D)|>C$.
On a alors $|\partial_ID|=|\partial(M\backslash D)|-|\partial M\backslash
\partial_ED|\geq C+|\partial_ED|-|\partial M|\geq C-\varepsilon|\partial M|$.
En choisissant $\varepsilon$ suffisamment petit par rapport à $C$ et 
$|\partial M|$, on a donc $|\partial_ID|\geq\varepsilon|\partial M|\geq
\varepsilon|\partial_ED|$ ce qui permet de conclure.

On montre maintenant la minoration de $\sigma_1$. Afin
d'alléger les notations on n'explicitera pas les références aux densités,
celles-ci étant généralement évidentes.
Soit $f$ une fonction propre sur $M$ de la première valeur propre 
$\sigma_1(M)$. On choisit son signe 
de sorte que $M^+=f^{-1}([0,+\infty[)$ vérifie $|M^+|\leq |M|/2$ (ou
$|\partial_EM^+|\leq |\partial M|/2$ si on suit la remarque~\ref{rem1}). 
La fonction $f$ restreinte à $M^+$ est alors la première fonction propre du 
problème de Steklov sur $M^+$ avec condition de Dirichlet sur $\partial_IM^+$, 
pour la valeur propre $\sigma_1$. Rappelons que problème de Steklov avec 
condition 
de Dirichlet sur une partie du bord de la variété est bien défini, que son 
spectre est strictement positif et qu'on a aussi la relation 
$\sigma_1=\int_{M^+}|\de f|^2/\int_{\partial M^+}f^2$ (cf. \cite{ag05}). 
Par conséquent on peut écrire, en utilisant 
l'inégalité de Cauchy-Schwarz :

\begin{eqnarray}
\sigma_1&=&\frac{\int_{M^+}|\de f|^2}{\int_{\partial M^+}f^2}
=\frac{\left(\int_{M^+}f^2\right)\left(\int_{M^+}|\de f|^2\right)}
{\left(\int_{M^+}f^2\right)\left(\int_{\partial M^+}f^2\right)}
 \geq\frac{\left(\int_{M^+}f|\de f|\right)^2}
{\left(\int_{M^+}f^2\right)\left(\int_{\partial M^+}f^2\right)}\nonumber\\
&\geq&\frac14\frac{\left(\int_{M^+}|\de(f^2)|\right)^2}
{\left(\int_{M^+}f^2\right)\left(\int_{\partial M^+}f^2\right)}
=\frac14\frac{\int_{M^+}|\de(f^2)|}{\int_{M^+}f^2}\cdot
\frac{\int_{M^+}|\de(f^2)|}{\int_{\partial M^+}f^2}.
\end{eqnarray}
En posant $D_t=f^{-1}([\sqrt t,+\infty[)$, la formule de la co-aire 
(\cite{kp08}, ch.~5) donne les trois relations
\begin{equation}
\int_{M^+}|\de(f^2)|=\int_{t\geq0}|\partial_ID_t|\de t,\ 
\int_{M^+}f^2=\int_{t\geq0}|D_t|\de t \textrm{ et }
\int_{\partial M^+}f^2=\int_{t\geq0}|\partial_ED_t|\de t.
\end{equation}
Comme $|D_t|\leq|M^+|\leq|M|/2$ (ou $|\partial_E D_t|\leq|\partial_E M^+|
\leq|\partial M|/2$) pour tout $t\geq0$, on en déduit que 
$\int_{M^+}|\de(f^2)|\geq h\int_{M^+}f^2$ et
$\int_{M^+}|\de(f^2)|\geq h'\int_{\partial M^+}f^2$,
et donc que $\sigma_1\geq \frac14 hh'$.
\end{proof}

\providecommand{\bysame}{\leavevmode ---\ }
\providecommand{\og}{``}
\providecommand{\fg}{''}
\providecommand{\smfandname}{\&}
\providecommand{\smfedsname}{\'eds.}
\providecommand{\smfedname}{\'ed.}
\providecommand{\smfmastersthesisname}{M\'emoire}
\providecommand{\smfphdthesisname}{Th\`ese}


\begin{thebibliography}{CESG11}

\bibitem[Agr05]{ag05}
{\scshape M.~S. Agranovich} -- {\og On a mixed {P}oincar\'e-{S}teklov {T}ype
  {S}pectral {P}roblem in a {L}ipschitz {D}omain\fg}, \emph{Russ. J. Math.
  Phys.} \textbf{13} (2005), no.~3, p.~239--244.

\bibitem[Ban80]{ba80}
{\scshape C.~Bandle} -- \emph{Isoperimetric inequalities and applications},
  Monographs and Studies in Mathematics, vol.~7, Pitman, 1980.

\bibitem[Bro85]{br85}
{\scshape R.~Brooks} -- {\og The bottom of the spectrum of a riemannian
  covering\fg}, \emph{J. Reine Angew. Math.} \textbf{357} (1985), p.~101--114.

\bibitem[Bus80]{bu80}
{\scshape P.~Buser} -- {\og On {C}heeger inequality $\lambda_1\geq h^2/4$\fg},
  in \emph{Geometry of the {L}aplace operator}, Proc. Sympos. Pure Math.,
  XXXVI, Amer. Math. Soc., 1980, p.~29--77.

\bibitem[CESG11]{cesg11}
{\scshape B.~Colbois, A.~El~Soufi {\normalfont \smfandname} A.~Girouard} --
  {\og Isoperimetric control of the {S}teklov spectrum\fg}, \emph{J. Funct.
  Anal.} \textbf{261} (2011), no.~5, p.~1384--1399.

\bibitem[Che70]{ch70}
{\scshape J.~Cheeger} -- {\og {A} lower bound for the smallest eigenvalue of
  the {L}aplacian\fg}, in \emph{{P}roblems in analysis ({P}apers dedicated to
  {S}alomon {B}ochner, 1969)}, Princeton Univ. Press, 1970.

\bibitem[CO97]{co97}
{\scshape S.-Y. Cheng {\normalfont \smfandname} K.~Oden} -- {\og Isoperimetric
  inequalities and the gap between the first and second eigenvalues of an
  {E}uclidean domain\fg}, \emph{J. Geom. Anal.} \textbf{7} (1997), no.~2,
  p.~217--239.

\bibitem[Esc97]{es97}
{\scshape J.~F. Escobar} -- {\og The geometry of the first non-zero {S}tekloff
  eigenvalue\fg}, \emph{J. Funct. Anal.} \textbf{150} (1997), no.~2,
  p.~544--556.

\bibitem[Esc99]{es99}
\bysame , {\og An isoperimetric {I}nequality and the first {S}teklov
  {E}igenvalue\fg}, \emph{J. Funct. Anal.} \textbf{165} (1999), no.~1,
  p.~101--116.

\bibitem[GP10]{gp10}
{\scshape A.~Girouard {\normalfont \smfandname} I.~Polterovich} -- {\og On the
  {H}ersch-{P}ayne-{S}chiffer inequalities for {S}teklov eigenvalues\fg},
  \emph{Functional Analysis and its Applications} \textbf{44} (2010), no.~2,
  p.~106--117.

\bibitem[Gu{\'e}04]{gu04}
{\scshape P.~Gu{\'e}rini} -- {\og Prescription du spectre du laplacien de
  {H}odge-de~{R}ham\fg}, \emph{Ann. scient. \'Ec. norm. sup. (4)} \textbf{37}
  (2004), no.~2, p.~270--303.

\bibitem[KP08]{kp08}
{\scshape S.~T. Krantz {\normalfont \smfandname} H.~R. Parks} --
  \emph{Geometric {I}ntegration {T}heory}, Birk\"auser, 2008.

\bibitem[SU90]{su90}
{\scshape J.~Sylvester {\normalfont \smfandname} G.~Uhlmann} -- {\og The
  {D}irichlet to {N}eumann map and applications\fg}, in \emph{Inverse problems
  in partial differential equations ({A}rcata, {CA}, 1989)}, SIAM, 1990,
  p.~101--139.

\bibitem[Uhl09]{uh09}
{\scshape G.~Uhlmann} -- {\og Electrical impedance tomography and
  {C}alder\'on's problem\fg}, \emph{Inverse Problems} \textbf{25} (2009),
  no.~12, p.~123011, 39.

\end{thebibliography}
\end{document}